# A CORRECTED QUANTITATIVE VERSION OF THE MORSE LEMMA

SÉBASTIEN GOUËZEL AND VLADIMIR SHCHUR

ABSTRACT. There is a gap in the proof of the main theorem in the article [Shc13a] on optimal bounds for the Morse lemma in Gromov-hyperbolic spaces. We correct this gap, showing that the main theorem of [Shc13a] is correct. We also describe a computer certification of this result.

## 1. INTRODUCTION

The Morse lemma is a fundamental result in the theory of Gromov-hyperbolic spaces. It asserts that, in a $\delta$-hyperbolic space, the Hausdorff distance between a $(\lambda, C)$-quasi-geodesic and a geodesic segment sharing the same endpoints is bounded by a constant $A(\lambda, C, \delta)$ depending only on $\lambda$, $C$ and $\delta$, and not on the length of the geodesic. Many proofs of this result have been given, with different expressions for $A$. An optimal value for $A$ (up to a multiplicative constant) has only been found recently in the article [Shc13a] by the second author, giving $A(\lambda, C, \delta) = K\lambda^2(C + \delta)$ for an explicit constant $K = 4(78 + 133/\log(2) \cdot \exp(157 \log(2)/28)) \sim 37723$.

Unfortunately, there is a gap in the proof of this theorem in [Shc13a], which was noticed by the first author while he was developing a library [Gou18] on Gromov-hyperbolic spaces in the computer assistant Isabelle/HOL. In such a process, all proofs are formalized on a computer, and checked starting from the most basic axioms. The degree of confidence reached after such a formal proof is orders of magnitude higher than what can be obtained by even the most diligent reader of referee, and indeed this process shed the light on the gap in [Shc13a]. The gap is on Page 829: the inequality $\sum_{i=1}^{n} e^{-X_i}(X_{i-1} - X_i) \leqslant \int_0^\infty e^{-t} \, dt$ goes in the wrong direction as the sequence $X_i$ is decreasing.

In this paper, we fix this gap. Here is the estimate we get.

**Theorem 1.1.** *Consider a $(\lambda, C)$-quasi-geodesic $Q$ in a $\delta$-hyperbolic space $X$, and $G$ a geodesic segment between its endpoints. Then the Hausdorff distance $HD(Q, G)$ between $Q$ and $G$ satisfies*
$$HD(Q, G) \leqslant 92\lambda^2(C + \delta).$$

Let us specify precisely the terms used in this statement, as there are small variations in the definitions in the literature. For us, a $(\lambda, C)$-quasi-geodesic is the image of a map $f$ from a compact interval to $X$ satisfying for all $x, y$ the inequalities
$$\lambda^{-1}|y - x| - C \leqslant d(f(x), f(y)) \leqslant \lambda|y - x| + C.$$

A map satisfying these inequalities is also called a $(\lambda, C)$-quasi-isometry. We also require $\lambda \geqslant 1$ and $C \geqslant 0$ in the definition. A geodesic segment is by definition a $(1, 0)$-quasi-geodesic.

---

We say that the space $X$ is $\delta$-hyperbolic if the Gromov product $(x,y)_w = (d(x,w)+d(y,w)-d(x,y))/2$ satisfies for all points $x, y, z, w$ the inequality

$$(x,z)_w \geqslant \min((x,y)_w, (y,z)_w) - \delta.$$

Finally, the Hausdorff distance $HD(Q, G)$ is the smallest number $r$ such that $G$ is included in the $r$-neighborhood of $Q$, and conversely.

The new proof of Theorem 1.1 has been completely formalized in Isabelle/HOL in [Gou18]. Therefore, the above theorem is certified. Here is this statement as proved in Isabelle/HOL.

```
theorem (in Gromov_hyperbolic_space) Morse_Gromov_theorem':
  fixes f::"real ⇒ 'a"
  assumes "lambda C-quasi_isometry_on {a..b} f"
          "geodesic_segment_between G (f a) (f b)"
  shows "hausdorff_distance (f`{a..b}) G ≤ 92 * lambda^2 * (C + deltaG(TYPE('a)))"
```

In this formal statement, `'a` is a type of class `Gromov_hyperbolic_space`. It corresponds to the space $X$ of Theorem 1.1, and the associated hyperbolicity constant is `deltaG(TYPE('a))`. Instead of talking of the quasi-geodesic $Q$, the formal statement is made in terms of its parametrization $f$, as the notion of endpoint of a quasi-geodesic is not really well defined. With this correspondence, the two statements directly correspond to each other.

Although the proof is more involved than the original argument in [Shc13a], the constant we get in the end is much better (92 instead of 37724). Indeed, we have tried to optimize the constant as much as we could, contrary to [Shc13a], keeping in mind the foundational nature of the library [Gou18]. This optimization owes a lot to the formalization process. It makes it possible to optimize locally one part of the proof, and see if it breaks other parts of the proofs by checking if the proof assistant complains that the proof is not correct any more, or if everything goes through. The certainty of the result also makes the optimization worth it, as we are sure not to have forgotten for example an edge case that would spoil the estimates.

Having a formalized certified proof raises interesting questions about the way to write mathematics. We do not need to convince a reader (or a referee!) that the result is correct, as we have already done the much more demanding task of convincing a computer, and the proof with all details can be read by the interested reader in [Gou18]. Rather, we have to convey the interesting ideas. We have decided to give all the precise statements we use (both in the traditional version and in the formalized version), but skip their proofs if they are small variations around results that are already available in the literature. For the main proof, we will explain (with as many details as in a traditional mathematical paper) a simplified version of the proof that gives the same statement as Theorem 1.1 but not caring much about the universal constants (this simplified argument gives the constant 2460 instead of 92 in Theorem 1.1). Then we will comment without entering in too many details on the various optimizations that can be done, leading to the above statement.

## 2. Proof of the main theorem

The proof uses the notion of quasiconvexity. We say that a subset $Y \subseteq X$ is $K$-quasiconvex if, for any $y_1, y_2 \in Y$, there exists a geodesic between $y_1$ and $y_2$ which is included in the $K$-neighborhood of $Y$. For instance, geodesics are 0-quasiconvex. The $r$-neighborhood of a 0-quasiconvex set is always $8\delta$-quasiconvex, see [CDP90, Proposition 10.1.2].

We follow the global strategy of [Shc13a] to prove Theorem 1.1, with a new more involved argument at a key technical step. Thanks to [BS00], we can assume without loss of generality that the space $X$ is geodesic. The quasi-geodesic $Q$ is by definition the image of a $(\lambda, C)$-quasi-isometric map $f : [u^-, u^+] \to X$. The statement for a general quasi-isometric map $f$ reduces to the one for a continuous quasi-isometric map $f$ thanks to the following approximation lemma, which is a version of [Shc13a, Lemma 9].

**Lemma 2.1.** *Any $(\lambda, C)$-quasi-isometry from a compact interval to a geodesic metric space can be uniformly approximated, up to $8C$, by a $(\lambda, 17C)$-quasi-geodesic with the same endpoints which is moreover $2\lambda$-Lipschitz.*[1]

Here is the formal version from which the above statement follows by taking $\Delta = 2$. Note that this formal statement requires that the endpoints of the original quasi-geodesic are separated by at least $2C$. This assumption is indeed necessary for the lemma to hold, but as it is satisfied in all applications it is most often overlooked in mathematical texts (see for instance [BH99, Lemma III.H.1.11] or [GdlH90, Lemma 5.8]).

```
lemma (in geodesic_space) quasi_geodesic_made_lipschitz:
  fixes c::"real ⇒ 'a"
  assumes "lambda C-quasi_isometry_on {a..b} c"
          "dist (c a) (c b) ≥ Delta * C" "Delta > 1"
  shows "∃d. continuous_on {a..b} d ∧ d a = c a ∧ d b = c b
          ∧ (∀x∈{a..b}. dist (c x) (d x) ≤ 4 * Delta * C)
          ∧ lambda ((8 * Delta + 1) * C)-quasi_isometry_on {a..b} d
          ∧ ((Delta/(Delta-1)) * lambda)-lipschitz_on {a..b} d"
```

From this point on, replacing $C$ by $17C$, we will assume that the $(\lambda, C)$-quasi-isometry $f$ is also continuous. Replacing the original hyperbolicity constant $\delta_0$ by a slightly larger constant $\delta$ (and letting $\delta$ tend to $\delta_0$ at the end of the argument), we can assume that the space is hyperbolic for a constant strictly smaller than $\delta$, and also that $\delta > 0$.

Consider $z \in [u^-, u^+]$. We want to estimate $d(f(z), G)$. We will prove an estimate of the form

$$(2.1) \qquad d(f(z), G) \leqslant K_0 + \frac{K_1}{K_2} \int_0^{u^+ - u^-} e^{-K_2 t} \, dt = K_0 + K_1 \cdot (1 - e^{-K_2(u^+ - u^-)}),$$

where $K_0$, $K_1$ and $K_2$ are suitable parameters that do not depend on $u^-$ and $u^+$. Both $K_0$ and $K_1$ will be of the form $K_i = k_i \lambda^2 (C+\delta)$, while $K_2$ will be of the form $K_2 = k_2/(\delta\lambda)$ where $k_0, k_1, k_2$ are explicit positive real constants. They will be defined in (2.4), (2.7) and (2.6). This estimate is proved inductively over the size of $u^+ - u^-$, reducing the estimate over $[u^-, u^+]$ to the estimate over a shorter interval $[v^-, v^+]$. We will have to show that the loss in this reduction process is controlled in terms of $K_1 e^{-K_2(v^+ - v^-)} - K_1 e^{-K_2(u^+ - u^-)}$, to conclude the proof of (2.1) by induction.

Let us first explain why this estimate concludes the proof. It implies that $d(f(z), G) \leqslant K_0 + K_1$. This proves that the image $Q$ of $f$ is included in the $(k_0 + k_1)\lambda^2(C+\delta)$-neighborhood of $G$. To get the estimate on the Hausdorff distance, one needs to show that $G$ is also included in a $k\lambda^2(C+\delta)$-neighborhood of $Q$ for some $k$. This follows from the previous estimate and a standard argument (see [BH99]) that we recall now. Consider a point $g \in G$. Denote by $Q^-$ the set of points on $Q$ that are within distance $(k_0+k_1)\lambda^2(C+\delta)$ of a point of $G$ in $[f(u^-),g]$, and by $Q^+$ the set of points on $Q$ that are within distance $(k_0+k_1)\lambda^2(C+\delta)$ of a point of $G$ in $[g, f(u^+)]$. The previous estimate implies that $Q = Q_1 \cup Q_2$. As $Q$ is connected, it follows that $Q_1 \cap Q_2 \neq \emptyset$. Denote by $f(z)$ a point in this intersection, and by $g^-$ and $g^+$ two points before and after $g$ on $G$, at distance at most $(k_0+k_1)\lambda^2(C+\delta)$ of $f(z)$. Using hyperbolicity in a triangle with vertices at $g^-, g^+, f(z)$ and the fact that $g$ is on a geodesic between $g^-$ and $g^+$, it follows that the distance between $g$ and $f(z)$ is at most $(k_0+k_1)\lambda^2(C+\delta)+\delta$. As $\lambda \geq 1$, this expression is bounded by $(k_0+k_1+1)\lambda^2(C+\delta)$. This concludes the argument, for the constant $k = k_0 + k_1 + 1$. We remind that [Shc13b] contains a stronger result (Theorem 3) claiming that the geodesic $G$ is included in an $A(\delta \log \lambda + C + \delta)$-neighborhood of the quasi-geodesic $Q$ with some universal constant $A$.

---

[1] provided that the distance between the endpoints of the initial quasi-geodesic is at least $2C$.

It remains to prove the estimate (2.1). The proof will use two parameters $L$ and $D$. For simplicity, let us take

(2.2) $$L = D = 100\delta.$$

We keep separate notations for $L$ and $D$ because we will want to optimize the choice of their values later.

*Case 1.* The case where $d(f(z), G) \leq L$ is trivial, as the estimate (2.1) holds if one takes $K_0$ large enough.

*Case 2.* Let us therefore assume $d(f(z), G) > L$. We will construct several points along $[u^-, z]$. To ease the reading, their order will correspond to the alphabetical order when possible.

Consider a projection $\pi_z$ of $f(z)$ on $G$, and $H$ a geodesic segment from $\pi_z$ to $f(z)$. Denote by $p : X \to H$ a closest-point projection on $H$. The idea is to project the quasi-geodesic $Q$ on $H$ and to consider the subpart $Q'$ of $Q$ that projects at distance at least $L$ of $\pi_z$. If one could show that $Q'$ is quantitatively shorter than $Q$ and that the distance from $f(z)$ to $\pi_z$ is controlled in terms of the distance from $f(z)$ to a geodesic joining the endpoints of $Q'$, then we would be in good shape to prove (2.1) inductively, deducing the estimate for $Q$ from the estimate for $Q'$. The real argument will be built around this naive idea, but in a more subtle way.

More precisely, consider two points $y^- \in [u^-, z]$ and $y^+ \in [z, u^+]$ such that the projections $p(f(y^-))$ and $p(f(y^+))$ are at distance roughly $L$ of $\pi_z$. In general, $p$ is not uniquely defined and not continuous, but this is almost the case up to $O(\delta)$ thanks to the hyperbolicity of the space. With the following standard lemma and recalling that $H$ is 0-quasiconvex as it is a geodesic, one can find $y^-$ and $y^+$ such that

(2.3) $$d(p(f(y^\pm)), \pi_z) \in [L - 4\delta, L].$$

**Lemma 2.2.** *A closest-point projection of a connected set on a $K$-quasiconvex subset of $X$ has gaps of size at most $4\delta + 2K$.*

Here is the formal version of this lemma we use.

```
lemma (in Gromov_hyperbolic_space_geodesic) quasi_convex_projection_small_gaps:
  assumes "continuous_on {a..(b::real)} f"
          "a ≤ b"
          "quasiconvex C G"
          "⋀t. t ∈ {a..b} ⟹ p t ∈ proj_set (f t) G"
          "delta > deltaG(TYPE('a))"
          "d ∈ {4 * delta + 2 * C..dist (p a) (p b)}"
  shows "∃t ∈ {a..b}. (dist (p a) (p t) ∈ {d - 4 * delta - 2 * C .. d})
                      ∧ (∀s ∈ {a..t}. dist (p a) (p s) ≤ d)"
```

Denote by $d^-$ (respectively $d^+$) the minimal distance of a point in $f([u^-, y^-])$ (respectively $f([y^+, u^+])$) to $H$. These distances are realized by two points $f(m^-)$ and $f(m^+)$.

*Case 2.1* Assume that $\max(d^-, d^+)$ is not large, say $\leqslant D + C$ where $D = 100\delta$ is the constant we have chosen in (2.2) and $C$ is the quasi-isometry parameter. This is again an easy case. Indeed, as the projections of $f(m^-)$ and $f(m^+)$ are within distance $L$ of $\pi_z$, one gets $d(f(m^-), f(m^+)) \leqslant 2D + 2C + L$. By quasi-isometry,

$$d(m^-, m^+) \leqslant \lambda(d(f(m^-), f(m^+)) + C) \leqslant \lambda(2D + 3C + L).$$

As $z$ is between $m^-$ and $m^+$, one gets in particular $d(m^-, z) \leqslant \lambda(2D + 3C + L)$. Then

$$d(f(z), \pi_z) \leqslant d(f(z), f(m^-)) + d(f(m^-), p(f(m^-))) + d(p(f(m^-)), \pi_z)$$
$$\leqslant (\lambda d(z, m^-) + C) + (D + C) + L \leqslant \lambda^2(3D + 5C + 2L).$$

This is compatible with the inequality (2.1) if one takes

(2.4) $$K_0 = 500\lambda^2(\delta + C).$$

*Case 2.2* Assume now that $\max(d^-, d^+) \geqslant D + C$, and $d^- \geqslant d^+$ for instance. This is the interesting case. The main step in the proof is the following lemma.

**Lemma 2.3.** *There exist two points $v \leqslant x$ in $[u^-, y^-]$ and a real number $d' \geqslant d^-$ such that*

(2.5) $$L - 74\delta \leqslant 4\sqrt{2}\lambda(x - v)e^{-d'\log(2)/(10\delta)}$$

*and $d(f(v), p(f(v))) \leqslant 4d'$.*

We do not state a formal version of this lemma as it only appears in the middle of the formal proof, not in lemma form. The numerology in the lemma (74 and $4\sqrt{2}$ and $\log(2)/10$ and 4) is of no importance: what only matters is that $L - 74\delta$ is positive, thanks to the choice of $L$ in (2.2), and that the other numbers are positive and fixed.

Let us show how to conclude the proof using the lemma. We have

$$m^+ - v = d(v, m^+) \leqslant \lambda(d(f(v), f(m^+)) + C)$$
$$\leqslant \lambda\Big(d(f(v), p(f(v))) + d(p(f(v)), p(f(m^+))) + d(p(f(m^+)), f(m^+)) + C\Big)$$
$$\leqslant \lambda(4d' + L + d^+ + C) \leqslant 6\lambda d',$$

as $L + C = D + C \leqslant d^- \leqslant d'$ and $d^+ \leqslant d^- \leqslant d'$. Therefore, taking

(2.6) $$K_2 = \log(2)/(60\delta\lambda),$$

the inequality (2.5) gives

$$L - 74\delta \leqslant 4\sqrt{2}\lambda(x-v)e^{-(m^+-v)\cdot\log(2)/(60\delta\lambda)} = \frac{4\sqrt{2}\lambda}{K_2} \cdot K_2(x-v)e^{-K_2(m^+-v)}$$

$$\leqslant \frac{4\sqrt{2}\lambda}{K_2}(e^{K_2(x-v)} - 1)e^{-K_2(m^+-v)} = \frac{4\sqrt{2}\lambda}{K_2}(e^{-K_2(m^+-x)} - e^{-K_2(m^+-v)})$$

$$\leqslant \frac{4\sqrt{2}\lambda}{K_2}(e^{-K_2(m^+-x)} - e^{-K_2(u^+-u^-)}).$$

Consider a new geodesic $G'$ between $f(x)$ and $f(m^+)$. Arguing by induction, we can assume that the estimate (2.1) has already been proved for $G'$, and we want to deduce it for $G$. Since both endpoints of $G'$ project within distance $L$ of $\pi_z$, one checks that the distance from $f(z)$ to $G$ is controlled by the distance from $f(z)$ to $G'$ (this is a version of [Shc13a, Lemma 5]). More specifically,

$$d(f(z), G) \leqslant d(f(z), G') + L + 4\delta.$$

Bounding $d(f(z), G')$ thanks to the induction assumption, and plugging in the estimate from the previous equation, we get

$$d(f(z), G) \leqslant K_0 + K_1(1 - e^{-K_2(m^+-x)}) + \frac{L+4\delta}{L-74\delta} \cdot \frac{4\sqrt{2}\lambda}{K_2}(e^{-K_2(m^+-x)} - e^{-K_2(u^+-u^-)}).$$

Let us take

(2.7) $$K_1 = \frac{L+4\delta}{L-74\delta} \cdot \frac{4\sqrt{2}\lambda}{K_2}.$$

Then the terms $K_1 e^{-K_2(m^+-x)}$ simplify in this equation, and we are left with

$$d(f(z), G) \leqslant K_0 + K_1(1 - e^{-K_2(u^+-u^-)}).$$

This is (2.1), as desired. This concludes the proof of Theorem 1.1. □

It remains to prove Lemma 2.3. The argument relies on the contracting properties of closest-point projections on quasiconvex sets. The first such basic statement is the following variation around [CDP90, Proposition 10.2.1].

**Lemma 2.4.** *Consider a $K$-quasiconvex subset $Y$ of $X$. Then projections $p_x$ and $p_y$ on $Y$ of two points $x$ and $y$ satisfy*

$$d(p_x, p_y) \leqslant \max(5\delta + 2K, d(x,y) - d(x,p_x) - d(y,p_y) + 10\delta + 4K).$$

This result expresses the classical fact that a geodesic from $x$ to $y$ essentially follows a geodesic from $x$ to $p_x$, then from $p_x$ to $p_y$, then from $p_y$ to $y$. The formal version follows.

```
lemma (in Gromov_hyperbolic_space_geodesic) proj_along_quasiconvex_contraction:
  assumes "quasiconvex C G" "px ∈ proj_set x G" "py ∈ proj_set y G"
  shows "dist px py ≤ max (5 * deltaG(TYPE('a)) + 2 * C)
          (dist x y - dist px x - dist py y + 10 * deltaG(TYPE('a)) + 4 * C)"
```

The second result we need is more sophisticated. Instead of a linear gain in terms of the distance to the set one projects on, as in the previous lemma, it gives an exponential gain in the upper bound, by a successive reduction process. It is proved by putting points along the path with gaps of size $10\delta$. Then, move by $5\delta$ towards $Y$: this reduces the distance between the points by $5\delta$ essentially thanks to the previous lemma. Then, discard half the points: this shows that by moving towards $Y$ by $5\delta$ the length of the path has been divided by 2. One can iterate this argument to get the exponential gain. We give a statement for the projection on quasiconvex sets as this is what we will need later on. This statement is proved in [Shc13a, Lemma 10] for the projection on a geodesic segment, but the case of a general quasiconvex set is analogous.

**Lemma 2.5.** *Consider a $(\lambda, C)$-quasi-geodesic path $f : [a, b] \to X$, everywhere at distance at least $D$ of a $K$-quasiconvex subset $Y$. Then, if $D$ is large enough, projections $p_a$ of $f(a)$ and $p_b$ of $f(b)$ on $Y$ satisfy the inequality*

$$d(p_a, p_b) \leqslant 2K + 8\delta + \max\left(5\delta, 4\sqrt{2}\lambda(b-a)\exp\left(-(D - K - C/2)\log(2)/(5\delta)\right)\right).$$

The formal version follows.

```
lemma (in Gromov_hyperbolic_space_geodesic) quasiconvex_projection_exp_contracting:
  assumes "quasiconvex K G"
          "⋀x y. x ∈ {a..b} ⟹ y ∈ {a..b}
                        ⟹ dist (f x) (f y) ≤ lambda * dist x y + C"
          "pa ∈ proj_set (f a) G" "pb ∈ proj_set (f b) G"
          "⋀t. t ∈ {a..b} ⟹ infdist (f t) G ≥ D"
          "D ≥ 15/2 * delta + K + C/2"
          "delta > deltaG(TYPE('a))"
          "C ≥ 0" "lambda ≥ 0" "a ≤ b"
  shows "dist pa pb ≤ 2*K + 8*delta + max (5*deltaG(TYPE('a)))
     ((4*exp(1/2 * ln 2)) * lambda * (b-a) * exp(-(D-K-C/2) * ln 2 / (5*delta)))"
```

Using these results, we can prove Lemma 2.3.

*Proof of Lemma 2.3.* For $k \geqslant 0$, let $V_k$ denote the $(2^k - 1)d^-$-neighborhood of $H$. These sets are all $8\delta$-quasiconvex. We recall that $p(f(x))$ is a projection of $f(x)$ on $H$. Let $p_k(x)$ denote the point on a fixed geodesic between $p(f(x))$ and $f(x)$ at distance $\min((2^k - 1)d^-, d(p(f(x)), f(x)))$ of $p(f(x))$. Then $p_k(x)$ is a projection of $f(x)$ on $V_k$, and moreover these projections are compatible in the following sense: for $k \leqslant \ell$, then $p_k(x)$ is a projection of $p_\ell(x)$ on $V_k$. Moreover, $p_0(x) = p(f(x))$.

We will do an inductive construction over $k$. This construction will have to stop at some step, where it will give the desired points. Until the argument stops, we will construct a point $x_k \in [u^-, y^-]$ such that

(2.8) $$d(p_k(u^-), p_k(x_k)) \geqslant L - 8\delta$$

and

(2.9) $$\text{for all } w \in [u^-, x_k], \, d(f(w), p_0(w)) \geqslant (2^{k+1} - 1)d^-.$$

Let us first check that this property holds for $k = 0$. Take $x_0 = y^-$. The point $\pi_z$ is a projection of $f(z)$ on the geodesic $G$ between $f(u^-)$ and $f(u^+)$. This does not imply that the projection $p_0(u^-)$ of $f(u^-)$ on the geodesic $H$ between $\pi_z$ and $f(z)$ is exactly at $\pi_z$ (contrary to the situation in the Euclidean plane), but by hyperbolicity one checks that $d(\pi_z, p_0(u^-)) \leqslant 4\delta$ (this is a version of [Shc13a, Lemma 3]). Since $d(\pi_z, p_0(y^-)) \in [L-4\delta, L]$ by (2.3) and $x_0 = y^-$, we deduce that $d(p_0(u^-), p_0(x_0)) \geqslant L - 8\delta$. This is (2.8). Moreover, by definition of $d^-$, the inequality (2.9) holds for $k = 0$.

Assume now that (2.8) and (2.9) hold at $k$. We will show that either we can find a pair of points that satisfy the conclusion of the lemma, or we can construct a point $x_{k+1}$ such that (2.8) and (2.9) hold at $k+1$.

As $V_k$ is $8\delta$-quasiconvex, we deduce from Lemma 2.2 that the gaps of the closest-point projection $p_k$ are bounded by $20\delta$. Therefore, we can find a point $x_{k+1} \in [u^-, x_k]$ whose projection on $V_k$ satisfies

$$(2.10) \qquad d(p_k(u^-), p_k(x_{k+1})) \in [22\delta, 42\delta],$$

and moreover all points $w \in [u^-, x_{k+1}]$ satisfy

$$(2.11) \qquad d(p_k(u^-), p_k(w)) \leqslant 42\delta.$$

There are two cases to consider:

*If there exists $v \in [u^-, x_{k+1}]$ with $d(f(v), p_0(v)) \leqslant (2^{k+2} - 1)d^-$.* Then we claim that the pair $(v, x_k)$ satisfies the conclusion of Lemma 2.3, for $d' = 2^k d^-$. First, the inequalities $d' \geqslant d^-$ and $d(f(v), p_0(v)) \leqslant 4d'$ hold by construction. Moreover, $d(p_k(v), p_k(x_k)) \geqslant L - 50\delta$ as $p_k(x_k)$ is far from $p_k(u^-)$ by (2.8), and $p_k(v)$ is close to $p_k(u^-)$ by (2.11). As all intermediate points are at distance at least $(2^{k+1} - 1)d^-$ of $V_0$ by (2.9), they are at distance at least $2^k d^-$ of $V_k$ and we can apply the exponential contraction lemma 2.5 with $D = 2^k d^-$. As $V_k$ is $8\delta$-quasiconvex, we get

$$L - 50\delta \leqslant d(p_k(v), p_k(x_k))$$
$$\leqslant 24\delta + \max\left(5\delta, 4\sqrt{2}\lambda(x_k - v)\exp\left(-(2^k d^- - 8\delta - C/2)\log(2)/(5\delta)\right)\right).$$

As $L - 50\delta > 29\delta$, the maximum has to be realized by the second term. Moreover, $2^k d^- - 8\delta - C/2 \geqslant (2^k d^-)/2 = d'/2$, as $d^- \geqslant D + C = 100\delta + C$. We obtain

$$(2.12) \qquad L - 74\delta \leqslant 4\sqrt{2}\lambda(x_k - v)\exp\left(-d'\log(2)/(10\delta)\right).$$

This concludes the proof in this case.

*Otherwise, $d(f(w), p_0(w)) \geqslant (2^{k+2} - 1)d^-$ for all $w \in [u^-, x_{k+1}]$.* In this case, (2.9) holds for $k+1$. Let us check that (2.8) also holds for $k+1$, by applying the projection lemma 2.4 to the points $p_{k+1}(u^-)$ and $p_{k+1}(x_{k+1})$, which project respectively to $p_k(u^-)$ and $p_k(x_{k+1})$ on $V_k$. As $V_k$ is $8\delta$-quasiconvex, this lemma gives

$$d(p_k(u^-), p_k(x_{k+1})) \leqslant \max(21\delta,$$
$$d(p_{k+1}(u^-), p_{k+1}(x_{k+1})) - d(p_{k+1}(u^-), p_k(u^-)) - d(p_{k+1}(x_{k+1}), p_k(x_{k+1})) + 42\delta).$$

As $d(p_k(u^-), p_k(x_{k+1})) \geqslant 22\delta$ by (2.10), the maximum has to be realized by the second term. Both distances $d(p_{k+1}(u^-), p_k(u^-))$ and $d(p_{k+1}(x_{k+1}), p_k(x_{k+1}))$ are equal to $2^k d^-$. We obtain
$$2 \cdot 2^k d^- - 20\delta \leqslant d(p_{k+1}(u^-), p_{k+1}(x_{k+1})).$$
As $d^- \geqslant D = 100\delta$, the left hand side is $\geqslant L - 8\delta = 92\delta$. This concludes the proof of (2.8), and of the induction.

Finally, if the conclusion of the lemma does not hold, then the induction will go on forever. Taking in particular $w = u^-$ in (2.9), we get $d(f(u^-), p_0(u^-)) \geqslant (2^{k+1} - 1) d^-$ for all $k$, a contradiction. □

Here are some ways to optimize the proof to get better constants. In addition to multiple minor optimizations, let us mention the main ones:

- The set $V_0$ is 0-quasiconvex, not only $8\delta$-quasiconvex. This means that estimates in the proof of Lemma 2.3 are better for $k = 0$. There is a different source of gain for $k > 0$, thanks to the factor $2^k$. Separating the two cases improves the final constant.
- There is an exponential gain in (2.12). One can spend some part of this gain, say $\exp(-(1-\alpha)d' \log(2)/(10\delta)) \leqslant \exp(-(1-\alpha)D\log(2)/(10\delta))$ to improve the multiplicative constant, and use the remaining part $\exp(-\alpha d' \log(2)/(10\delta))$ for the induction (for a suitable value of $\alpha$).
- Instead of formulating the induction in terms of the distance from $f(z)$ to a geodesic $G$ between $f(u^-)$ and $f(u^+)$, it is more efficient to induce over the Gromov product $(f(u^-), f(u^+))_{f(z)}$ (which coincides with the distance $d(f(z), G)$ up to $2\delta$) as most inequalities are done in terms of Gromov products. The main interest of this change is that, with the current argument, the point $f(u^-)$ projects on $H$ between $\pi_z$ and $f(z)$ within distance $4\delta$ of $\pi_z$, which means there is a small loss. With the Gromov product approach, let $m$ denote the point on $G$ which is opposite to $f(z)$ in the triangle $[f(z), f(u^-), f(u^+)]$, i.e., it is on $G$ at distance $(f(z), f(u^+))_{f(u^-)}$ of $f(u^-)$ and at distance $(f(z), f(u^-))_{f(u^+)}$ of $f(u^+)$. Let $\pi_z$ denote the point on a geodesic $H$ from $f(z)$ to $m$ at distance $(f(u^-), f(u^+))_{f(z)}$ of $f(z)$. This point is within distance $2\delta$ of $m$. It turns out that the projection of $f(u^-)$ on $H$ is between $m$ and $\pi_z$, i.e., opposite from $f(z)$. The above loss is suppressed in this approach.
- Finally, one can choose freely $L$, $D$ and $\alpha$ within some range. In particular, $L$ and $D$ do not have to coincide. One can optimize numerically over these parameters to get the best possible bound. In the end, we take $L = 18\delta$ and $D = 55\delta$ and $\alpha = 12/100$ to get the value 92 in Theorem 1.1.

Laboratoire Jean Leray, CNRS UMR 6629, Université de Nantes, 2 rue de la Houssinière, 44322 Nantes, France
  *Email address*: sebastien.gouezel@univ-nantes.fr

Departments of Integrative Biology and Statistics, University of California Berkeley 4098 Valley Life Sciences Building (VLSB) Berkeley, CA 94720-3140, USA
  *Email address*: vlshchur@gmail.com